\documentclass[12pt]{article}
\usepackage[latin1]{inputenc}
\usepackage{amsmath}
\usepackage{amsfonts}
\usepackage{amssymb}
\usepackage{amsthm}

\usepackage{titlesec}
\titleformat{\subsection}[hang]{\normalfont\bfseries}{\thesubsection}{1em}{}

\usepackage{titlesec}

\titlespacing\section{0pt}{3.5ex plus 0.5ex minus .2ex}{0.3ex plus .2ex}
\titlespacing\subsection{0pt}{2.5ex plus 0.5ex minus .2ex}{0.3ex plus .2ex}
\titlespacing\subsubsection{0pt}{2.5ex plus 0.5ex minus .2ex}{0.3ex plus .2ex}

\usepackage{mathrsfs} 

\usepackage[alphabetic,lite]{amsrefs} 

\usepackage{fullpage}
\usepackage{setspace}
\usepackage{hyperref}
\usepackage{color}
\usepackage{enumitem}
\usepackage{ulem} 

\usepackage{marginnote}

\usepackage[all,cmtip]{xy} 

\usepackage{fancyhdr}
\newtheorem{Thm}{Theorem}[section]

\newtheorem{Lemma}[Thm]{Lemma}

\theoremstyle{definition}
\newtheorem{Def}[Thm]{Definition}
\newtheorem{Rem}[Thm]{Remark}

\newcommand{\bC}{\mathbb{C}}

\newcommand{\bF}{\mathbb{F}}

\newcommand{\bQ}{\mathbb{Q}}
\newcommand{\bR}{\mathbb{R}}

\newcommand{\bZ}{\mathbb{Z}}

\newcommand{\cP}{\mathcal{P}}

\newcommand{\fg}{\mathfrak{g}}

\newcommand{\fs}{\mathfrak{s}}
\newcommand{\ft}{\mathfrak{t}}

\newcommand{\sA}{\mathscr{A}}
\newcommand{\sB}{\mathscr{B}}

\DeclareMathAlphabet{\mathpzc}{OT1}{pzc}{m}{it}

\newcommand{\Proof}{\textbf{Proof.\\}}

\newcommand{\ra}{\rightarrow}

\newcommand{\wt}{\widetilde}

\newcommand{\ov}{\overline}

\DeclareMathOperator{\Hom}{Hom}			
\DeclareMathOperator{\Id}{Id}		  	
\DeclareMathOperator{\GL}{GL}		  	
\DeclareMathOperator{\End}{End}			
\DeclareMathOperator{\Ind}{Ind}			
\DeclareMathOperator{\Sp}{Sp}		  	
\DeclareMathOperator{\Lie}{Lie}  	  
\DeclareMathOperator{\ind}{ind}  	  
\DeclareMathOperator{\cind}{c-ind}  	  
\DeclareMathOperator{\val}{val}  	  

\newcommand{\mmu}{\boldsymbol\mu}

\AtEndDocument{\bigskip{\footnotesize%
		\textsc{University of Cambridge, Cambridge, UK and Duke University, Durham, NC, USA} \par  
		\textit{Mailing address}: Trinity College, Cambridge, CB2 1TQ, UK \par
		\textit{E-mail address}: \texttt{fintzen@maths.cam.ac.uk} and \texttt{fintzen@math.duke.edu}
}}

\newcommand{\repKYu}{\wt \rho}

\newcommand{\KYu}{\wt K}

\newcommand{\varphil}{\varphi_\ell}
\newcommand{\barFell}{\overline\bF_\ell}
\newcommand{\Zlplus}{\overline\bZ_\ell^+}

\begin{document}
\author{Jessica Fintzen}
\title{\textcolor{white}{.}\\[-2cm] Tame cuspidal representations in non-defining characteristics} 
\date{}
\maketitle

\begin{abstract}
  Let $F$ be a non-archimedean local field of residual characteristic $p \neq 2$. Let $G$ be a (connected) reductive group that splits over a tamely ramified field extension of $F$. We show that a construction analogous to Yu's construction of complex supercuspidal representations yields smooth, irreducible, cuspidal representations over an arbitrary algebraically closed field $R$ of characteristic different from $p$. Moreover, we prove that this construction provides all smooth, irreducible, cuspidal $R$-representations if $p$ does not divide the order of the Weyl group of $G$. 
  \end{abstract}

{
	\renewcommand{\thefootnote}{}  
	\footnotetext{MSC2010: 22E50; 20C20, 20G05, 20G25} 
	\footnotetext{Keywords: representations of reductive groups over non-archimedean local fields, 
		 $p$-adic groups, cuspidal  $\ell$-modular representations}
	\footnotetext{The author was partially supported by NSF Grant DMS-1802234 / DMS-2055230 and a Royal Society University Research Fellowship.}
}

\setcounter{tocdepth}{1}

\tableofcontents

\newpage

\section{Introduction}
 Let $F$ be a non-archimedean local field of residual characteristic $p \neq 2$, and let $G$ be a (connected) reductive group that splits over a tamely ramified field extension of $F$. If $p$ does not divide the order of the (absolute) Weyl group of $G$, then we have an exhaustive construction of all smooth complex irreducible supercuspidal representations  of the $p$-adic group $G(F)$ (\cite{Yu, Fi-Yu-works, Kim, Fi-exhaustion}), which has been widely used to date.
 For number theoretic applications, it is crucial to also have a construction of smooth representations with coefficients in an algebraically closed field $R$ of positive characteristic $\ell \neq p$. 
   In this paper, we explain how to construct smooth $R$-representations using an analogous construction to the complex case. We prove that the resulting representations are irreducible and cuspidal, and we show that if $p$ does not divide the order of the Weyl group of $G$, then this construction yields all smooth, irreducible, cuspidal  $R$-representations  of $G(F)$.
 
While this paper is only concerned with algebraically closed fields $R$ as coefficient fields of our representations, Henniart and Vigneras (\cite{Henniart-Vigneras}) together with forthcoming work of Deseine use our results to obtain a similar result for non-algebraically closed coefficient fields. 
 
Modular representation theory has recently received new attention, e.g. also for the study of the local Langlands correspondence in families and global applications. However, $\ell$-modular representations have so far only been constructed for $\GL_n$ by Vignéras (\cite{Vigneras-book}) in 1996 following Bushnell--Kutzko's construction of complex representations (\cite{BK}), and recently for classical groups by Kurinczuk and Stevens (\cite{Kurinczuk-Stevens}) based on Stevens' earlier construction of complex representations (\cite{Stevens}).
 
 Our approach is closely following the setting with complex coefficients based on Yu's construction (\cite{Yu}). However, not all arguments for complex representations work also for mod-$\ell$ representations; there are some important differences, e.g. the $R$-representations of finite reductive groups over $\bF_p$ are no longer necessarily completely reducible, and in order to show that the compactly induced representation $\cind_{\KYu}^{G(F)}\repKYu$ is irreducible for an irreducible $R$-representation $\repKYu$, it does not always suffice to only show that every element that intertwines $\repKYu$ is contained in $\KYu$.
 Moreover, the cuspidal support does not in general decompose the category of $R$-representations into blocks (analogous to the Bernstein blocks in the complex setting). Hence the theory of types, which we used in the proof of exhaustion of complex supercuspidal representations, cannot be applied to this setting.
 Therefore we are required to provide some additional arguments in the mod-$\ell$ setting, which is the focus of this paper.
 In order to keep this paper short, we will not repeat all the arguments from the setting with complex coefficients that carry over to the mod-$\ell$ setting, but instead we provide precise references for the reader to read the original arguments that work equally well with $R$-coefficients and we focus on the additional arguments that need to be added.

	\textbf{Conventions and notation.} 
		All reductive groups in this paper are required to be connected.
		
		Let $F$ be a non-archimedean local field of residual characteristic $p \neq 2$ with valuation map $\val: F^* \ra \bZ $.  We denote by $\cP$ the maximal ideal of the ring of integers of $F$. 
			For a reductive group $G$ defined over $F$ we denote by $\sB(G,F)$ the (enlarged) Bruhat--Tits building (\cites{BT1,BT2}) of $G$ over $F$, by $Z(G)$ the center of $G$ and by $G^{\mathrm{der}}$ the derived subgroup of $G$. If $T$ is a maximal torus of $G$, we write $\Phi(G,T)$ for the roots of $G_{\ov F}$ with respect to $T_{\ov F}$, where $\ov F$ denotes a separable closure of $F$.  We let $\wt \bR=\bR \cup \{ {r+} \, | \, r \in \bR\}$ with its usual order, i.e. for $r$ and $s$ in $\bR$ with $r<s$, we have $r<r+<s<s+$. For $r \in \wt \bR_{\geq 0}$, we write $G_{x,r}$ for the Moy--Prasad filtration subgroup of $G(F)$ of depth $r$ at a point $x \in \sB(G,F)$. For $r \in \wt \bR$, we write $\fg_{x,r}$ or $\Lie(G)_{x,r}$ for the Moy--Prasad filtration submodule of  $\fg=\Lie G(F)$ of depth $r$ at $x$, and $\fg^*_{x,r}$ or $\Lie^*(G)_{x,r}$ for the Moy--Prasad filtration submodule of depth $r$ at $x$ of the linear dual $\fg^*$ of $\fg$.  If $x \in \sB(G,F)$, then we denote by $[x]$ its image in the reduced Bruhat--Tits building and we write 
	$G_{[x]}$ for the stabilizer of $[x]$ in $G(F)$. 
	
	We call a subgroup $G'$ of $G$ (defined over $F$) a twisted Levi subgroup of $G$ if $(G')_E=G'\times_F E'$ is a Levi subgroup of $G_E$ for some (finite) field extension $E$ of $F$. If $G'$ splits over a tamely ramified field extension of $F$, then, using (tame) Galois descent, we obtain an embedding of the corresponding Bruhat--Tits buildings $\sB(G',F) \hookrightarrow \sB(G,F)$. This embedding is only unique up to some translation, but its image is unique, and we will identify $\sB(G',F)$ with its image in $\sB(G,F)$. All constructions in this paper are independent of the choice of such an identification.

	If $K$ is a subgroup of $G$, $g \in G$, and $\rho$ a representation of $K$, then we write $^gK$ to denote $gKg^{-1}$ and define ${^g\rho}(x)=\rho(g^{-1}xg)$ for $x \in K \cap {^gK}$. 

Throughout the paper we let $\ell$ be a prime number different from $p$ and let $R$ be an algebraically closed field of characteristic $\ell$. We fix an additive character $\varphil: F \ra R^*$ of $F$ of conductor $\cP$ and a reductive group $G$ that is defined over our non-archimedean local field $F$ and that splits over a tamely ramified field extension of $F$. All representations of $G(F)$ in this paper have $R$-coefficients unless specified otherwise and are required to be smooth.

\textbf{Acknowledgements. } 
The author thanks Jeffrey Adler, Guy Henniart, Tasho Kaletha, Maarten Solleveld, Loren Spice and Marie-France Vigneras for discussions related to this paper. The author is also grateful for inspiration from the work of Robert Kurinczuk and Shaun Stevens for classical groups.

\section{Constructions of representations} \label{Section-mod-ell}
 In this section we explain how the construction of representations of supercuspidal complex representations by J.-K. Yu (\cite{Yu}) can be adapted to construct representations of $G(F)$ with coefficients in $R$. In Section \ref{Section-cuspidal} we will prove that the resulting representations are irreducible and cuspidal. 

\subsection{Generic characters}
We first introduce generic characters based on Yu's definition (\cite[\S8 and \S9]{Yu}) adapted to the mod-$\ell$ setting as this also allows us to address a small unclarity pointed out in \cite[Remark~4.1.3]{FKS}.
Let $G' \subsetneq G$ be a twisted Levi subgroup that splits over a tame extension, and denote by $(\Lie^*(G'))^{G'}$ the subscheme of $\Lie^*(G')$ fixed by (the dual of) the adjoint action of $G'$. 
\begin{Def}
	An elemente $X$ of $(\Lie^*(G'))^{G'}(F) \subset \Lie^*(G')(F)$ is called \textit{$G$-generic of depth $r \in \bR$} (or \textit{$(G, G')$-generic of depth $r$}) if the following three conditions hold.
	\begin{itemize}
		\item[\textbf{(GE0)}] For some (equivalently, every, see Lemma \ref{Lemma-for-all-x}) point $x \in \sB(G', F)$, we have $X \in \Lie^*(G')_{x,-r}$.
		\item[\textbf{(GE1)}] \textbf{(GE1)} of \cite[\S8]{Yu} holds, i.e. $\val(X(H_\alpha))=-r $ for all $\alpha \in \Phi(G, T) \setminus \Phi(G',T)$ for some maximal torus $T$ of $G'$, where $H_\alpha=\Lie(\alpha^\vee)(1)$ with $\alpha^\vee$ the coroot of $\alpha$. 
		\item[\textbf{(GE2)}] \textbf{(GE2)} of \cite[\S8]{Yu} holds, where we refer the reader to \cite{Yu} for details. Note that by \cite[Lemma~8.1]{Yu} \textbf{(GE1)} implies \textbf{(GE2)} if $p$ is not a torsion prime for the dual root datum of $G$.
	\end{itemize}
\end{Def}

\begin{Rem}
	Instead of condition \textbf{(GE0)}, J.-K. Yu identifies $(\Lie^*(G'))^{G'}$ with $\Lie^*(Z(G'))^\circ$ and requires $X$ to lie in $\Lie^*(Z(G'))^\circ_{-r}$. However, this identification does not always hold true, see \cite[Remark~4.1.3]{FKS}. 
\end{Rem}

\begin{Lemma} \label{Lemma-for-all-x}
	Let $X \in (\Lie^*(G'))^{G'}(F)$ and suppose that there exists a point $x \in \sB(G', F)$ such that $X \in \Lie^*(G')_{x,-r}$. Then $X \in \Lie^*(G')_{y,-r}$ for all points $y \in \sB(G', F)$. \\ 
	In particular, if $X$ is $(G, G')$-generic of depth $r$, then $X \in \Lie^*(G')_{y,-r}$ for all points $y \in \sB(G', F)$.
\end{Lemma}
\Proof
Without loss of generality (by replacing $F$ by a finite tame extension) we may assume that $G'$ is split. Let $y \in \sB(G', F)$, let $T$ be a maximal split torus of $G'$ whose apartment $\sA(T,F)$ contains $y$, and let $g \in G'(F)$ such that $g.x \in \sA(T,F)$. Since $X$ is fixed by the action of $g$, we have $X \in \Lie^*(G')_{g.x,-r}$. Moreover, we have a decomposition $\fg=\ft \oplus \ft^\perp$, where $\ft^\perp$ is the sum over the root subspaces of $\fg$ with respect to $T$. Since $X \in (\Lie^*(G'))^{G'}(F)$, the restriction of $X$ to $\ft^\perp$ is trivial, and therefore $X \in \Lie^*(G')_{g.x,-r}$ implies that $X \in \Lie^*(G')_{y,-r}$ as $g.x$ and $y$ are both in the apartment of $\sA(T,F)$.
\qed

Following \cite[\S9]{Yu} we define generic characters as follows.

\begin{Def} Let $x \in \sB(G',F)$ and $r \in \bR_{> 0}$.
	A character $\phi$ of $G'(F)$ is called \textit{$G$-generic (or $(G, G')$-generic) relative to $x$ of depth $r$} if $\phi$ is trivial on $G'_{x,r+}$, non-trivial on $G'_{x,r}$ and the restriction of $\phi$ to $G'_{x,r}/G'_{x,r+}\simeq \fg'_{x,r}/\fg'_{x,r+}$ is given by $\varphil \circ X$ for some $(G, G')$-generic element 
	$X$ of depth $r$.
\end{Def}

\subsection{The input for the construction of cuspidal representations} \label{Section-input-mod-ell}
The input for our construction of $R$-representations is the mod-$\ell$ analogue of the input for the constructions of supercuspidal complex representations, where we follow the set-up from Section 2.1 of \cite{Fi-Yu-works}. More precisely, the input is a tuple $((G_i)_{1 \leq i \leq n+1}, x, (r_i)_{1 \leq i \leq n}, \rho, (\phi_i)_{1 \leq i \leq n})$ for some non-negative integer $n$ where
\begin{enumerate}[label=(\alph*),ref=\alph*]
	\item $G=G_1 \supseteq G_2 \supsetneq G_3 \supsetneq \hdots \supsetneq G_{n+1}$ are twisted Levi subgroups of $G$ that split over a tamely ramified extension of $F$
	\item  $x \in \sB(G_{n+1},F)\subset \sB(G,F)$ 
	\item $r_1 > r_2 > \hdots > r_n >0$ are real numbers
	\item \label{item-rho-ell} $\rho$ is an irreducible $R$-representation of $(G_{n+1})_{[x]}$ that is trivial on $(G_{n+1})_{x,0+}$ 
	\item $\phi_i$, for $1 \leq i \leq n$, is an $R$-valued character of $G_{i+1}(F)$ of depth $r_i$ that is trivial on $(G_{i+1})_{x,r_i+}$ 
\end{enumerate}
satisfying the following conditions 
\begin{enumerate}[label=(\roman*),ref=\roman*]
	\item  $Z(G_{n+1})/Z(G)$ is anisotropic
	\item the image of the point $x$ in $\sB(G_{n+1}^{\mathrm{der}},F)$ is a vertex
	\item $\rho|_{(G_{n+1})_{x,0}}$ is a cuspidal representation of $(G_{n+1})_{x,0}/(G_{n+1})_{x,0+}$
	\item $\phi_i$ is $(G_i, G_{i+1})$-generic relative to $x$  of depth $r_i$ for all $1 \leq i \leq n$ with $G_i \neq G_{i+1}$
\end{enumerate}

\subsection{The construction of $R$-representations} \label{Section-construction-mod-ell}
We begin with a lemma that allows us to define a Weil--Heisenberg representation with $\barFell$-coefficients, which we then view as an $R$-representation via base change along an embedding $\barFell \hookrightarrow R$. In order to state the lemma, we denote by $\overline \bQ_\ell$ an algebraic closure of the $l$-adic numbers, by $\overline \bZ_\ell$ the integral closure of the algebraic integers $\bZ_\ell$ in $\overline \bQ_\ell$, and we write $\Zlplus$ for the maximal ideal of $\overline \bZ_\ell$. We let $\mmu_p$ be the subgroup of $\overline \bF_\ell^*$ consisting of the $p$-th roots of unity in  $\overline \bF_\ell^*$, and we use the Teichmüller lift to also view $\mmu_p$ as a subgroup of $\overline \bZ_\ell^* \subset \overline \bQ_\ell^*$.

\begin{Lemma}\label{Lemma-Weil-Heisenberg-mod-ell}
	Let $V$ be a symplectic $\bF_p$-vector space, let $V^\sharp=V \ltimes \bF_p$ be the corresponding Heisenberg $p$-group and $\varphi$ a non-trivial character of the center of $V^\sharp$ with values in $\mmu_p \subset \overline \bQ_\ell^*$. Let $(\omega, V_{\omega})$ denote a corresponding Weil--Heisenberg representation of $\Sp(V) \ltimes V^\sharp$ with coefficients in $\overline\bQ_\ell$ (obtained from \cite{Gerardin} via a fixed isomorphism $\bC \simeq \overline\bQ_\ell$). Then $V_{\omega}$ admits an  $\Sp(V) \ltimes V^\sharp$-stable $\overline\bZ_\ell$-lattice $L_{\omega}$, 
	 and the isomorphism class of the resulting $\barFell$-representation $\overline V_{\omega}:=L_{\omega}/(L_{\omega}\otimes_{\ov\bZ_\ell}\Zlplus)$ of $\Sp(V) \ltimes V^\sharp$ does not depend on the choice of $L_{\omega}$.
	
	Moreover, the representation $\overline V_{\omega}$ restricted to $V^\sharp$ is irreducible.
\end{Lemma}
\Proof As $\Sp(V) \ltimes V^\sharp$  is a finite group, the representation $V_{\omega}$ admits an  $\Sp(V) \ltimes V^\sharp$-stable $\overline\bZ_\ell$-lattice $L_{\omega}$ (see e.g. \cite[I.9.4]{Vigneras-book}). Since $V^\sharp$ is a $p$-group and $p \neq \ell$, the representation $\overline V_{\omega}$ restricted to $V^\sharp$ is irreducible. Hence also the representation $\overline V_{\omega}$  of $\Sp(V) \ltimes V^\sharp$ is irreducible and therefore uniquely (up to isomorphism) determined by its Brauer character. Since the Brauer character of $\overline V_{\omega}$ is the restriction of the character of $V_{\omega}$ to the $\ell$-regular elements of $\Sp(V) \ltimes V^\sharp$, the isomorphism class of the representation  $\overline V_{\omega}$ does not depend on the choice of $L_{\omega}$. \qed

We call the resulting $\barFell$-representation (and $R$-representation) from Lemma \ref{Lemma-Weil-Heisenberg-mod-ell} the \textit{mod-$\ell$ Weil--Heisenberg representation}. Using the mod-$\ell$ Weil--Heisenberg representation instead of the complex Weil--Heisenberg representation and our input from Section \ref{Section-input-mod-ell}, we can now imitate Yu's construction of supercuspidal representations (\cite{Yu}), see also Section 2.2 of \cite{Fi-Yu-works}, to obtain a compact-mod-center open subgroup $\KYu$ of $G(F)$ and a smooth representation $\repKYu$ of $\KYu$ such that the compact induction $\cind_{\KYu}^{G(F)}\repKYu$ yields a smooth $R$-representation of $G(F)$. Here
 \begin{eqnarray*}
	\KYu& =&(G_1)_{x,\frac{r_1}{2}}(G_2)_{x,\frac{r_2}{2}} \hdots (G_n)_{x,\frac{r_n}{2}}(G_{n+1})_{[x]} 
\end{eqnarray*}
and the representation $\repKYu$ of $\KYu$ is given by a tensor product $\rho \otimes \kappa$, where $\rho$ also denotes the extension of $\rho$ from $(G_{n+1})_{[x]}$ to $\KYu$ that is trivial on $(G_1)_{x,\frac{r_1}{2}}(G_2)_{x,\frac{r_2}{2}} \hdots (G_n)_{x,\frac{r_n}{2}}$, and $\kappa$ is a representation constructed from the $\phi_i$ using the mod-$\ell$ Weil--Heisenberg representation instead of the complex version, see \cite[Section~2.2]{Fi-Yu-works} for details.

\section{Irreducibility and cuspidality} \label{Section-cuspidal}

\begin{Thm} \label{Thm-mod-ell}
	The smooth $R$-representation $\cind_{\KYu}^{G(F)}\repKYu$ is irreducible and cuspidal.
\end{Thm}
We first show that a result by Gérardin (\cite{Gerardin}) for the complex Weil--Heisenberg representations also holds for the mod-$\ell$ Weil--Heisenberg representations.

\begin{Lemma} \label{Lemma-Gerardin-mod-ell}
	Let $V$ be a symplectic vector space over $\bF_p$, and let $V^+$ be a totally isotropic subspace of $V$ with orthogonal complement $V^+ \oplus V^0 \subset V$. We write $P$ for the (maximal) parabolic subgroup of $\Sp(V)$ that preserves the subspace $V^+$ (and hence also $V^+ \oplus V^0$) and obtain a surjection $P \twoheadrightarrow \Sp(V_0)$ by composing restriction to $V_0$ with the projection from $V^+ \oplus V^0$ to $V^0$ with kernel $V^+$.  
	Let $\varphi$ be a non-trivial character of $\bF_p$ with values in $\mmu_p \subset \overline\bF_\ell^*$, where we also view  $  \mmu_p \subset \overline \bQ_\ell^*$ via the Teichmüller lift. We denote by $\overline V_{\omega}$ and $\overline{V^0}_{\omega}$ the corresponding mod-$\ell$ Weil--Heisenberg representations of  $\Sp(V) \ltimes V^\sharp$ and  $\Sp(V^0) \ltimes (V^0)^\sharp$, respectively. Using the above projection from $P$ to $\Sp(V^0)$ and the projection from $V^+ \times (V^0)^\sharp$ to $(V^0)^\sharp$ with kernel $V^+$, we obtain a resulting action of $P \ltimes (V^+ \times (V^0)^\sharp)$ on $\overline{V^0}_{\omega}$.  Then
	$$\overline V_{\omega}|_{P \ltimes V^\sharp} \simeq \Ind_{P \ltimes (V^+ \times (V^0)^\sharp)}^{P \ltimes V^\sharp} \left(\overline{V^0}_{\omega} \otimes ({\barFell}_{\chi^{V^+}} \ltimes 1)\right), $$
	where ${\barFell}_{\chi^{V^+}}$ is the one dimensional vector space $\barFell$ on which the action of $p \in P$ is given by $\chi^{V^+}(p) = \det(p|_{V^+})^{(p-1)/2}$ and 1 denotes the trivial representation of $V^+ \times (V^0)^\sharp$.
\end{Lemma}
\Proof
By \cite[Theorem~2.4.(b)]{Gerardin}  	\footnote{The statement of \cite[Theorem~2.4.(b)]{Gerardin} omits the factor $\chi^{V_i^+} \ltimes 1$ (denoted by $\chi^{E_+}$ in \cite[Theorem~2.4.(b)]{Gerardin}), which is a typo.}
the restriction of the complex Weil--Heisenberg representation $V_{\omega}$ attached to $\varphi$ from $\Sp(V) \ltimes V^\sharp$ to $P \ltimes V^\sharp$ is 
given by 
$\Ind_{P \ltimes ( V^+ \times (V^0)^\sharp)}^{P \ltimes V^\sharp}\left( V_{\omega}^0 \otimes (\bC_{\chi^{V^+}} \ltimes 1)\right)$.
 We now view the $\bC$-representations as $\overline\bQ_\ell$-representations via a fixed isomorphism $\bC \simeq \overline\bQ_\ell$. Using an  $\Sp(V^0) \ltimes (V^0)^\sharp$-stable $\overline\bZ_\ell$-lattice $L_{\omega}^0$ in the $\overline\bQ_\ell$-Weil--Heisenberg representation $V_{\omega}^0$, we obtain a $P \ltimes V^\sharp$-stable $\overline\bZ_\ell$-lattice of $\Ind_{P \ltimes  ( V^+ \times (V^0)^\sharp)}^{P \ltimes V^\sharp} V_{\omega}^0 \otimes ({\overline\bQ_\ell}_{\chi^{V_1^+}} \ltimes 1)$ so that reduction mod~$\ell$ yields
$\Ind_{P \ltimes  ( V^+ \times (V^0)^\sharp)}^{P \ltimes V^\sharp} \left(\overline{V^0}_{\omega} \otimes ({\barFell}_{\chi^{V_1^+}} \ltimes 1)\right)$. Note that $\Ind_{P \ltimes  ( V^+ \times (V^0)^\sharp)}^{P \ltimes V^\sharp} \left(V_{\omega}^0 \otimes ({\overline\bQ_\ell}_{\chi^{V_1^+}} \ltimes 1)\right)$ restricted to $V^\sharp$ is isomorphic to the irreducible Heisenberg-representation $V_{\omega}$ of $V^\sharp$. Since $V^\sharp$ is a $p$-group with $p \neq \ell$, we deduce that $\Ind_{P \ltimes  ( V^+ \times (V^0)^\sharp)}^{P \ltimes V^\sharp} \left(\overline{V^0}_{\omega} \otimes ({\barFell}_{\chi^{V_1^+}} \ltimes 1)\right)$ is an irreducible representation of $P \ltimes V^\sharp$ and is therefore determined (up to isomorphism) by its Brauer character. The Brauer character of $\Ind_{P \ltimes  ( V^+ \times (V^0)^\sharp)}^{P \ltimes V^\sharp} \left(\overline{V^0}_{\omega} \otimes ({\barFell}_{\chi^{V_1^+}} \ltimes 1)\right)$  is the restriction of the character of $\Ind_{P \ltimes  ( V^+ \times (V^0)^\sharp)}^{P \ltimes V^\sharp} V_{\omega}^0 \otimes ({\overline\bQ_\ell}_{\chi^{V_1^+}} \ltimes 1)$, and therefore also the restriction of the  character of  the Weil--Heisenberg representation $V_{\omega}$ to the $\ell$-regular elements of $P \ltimes V^\sharp$. 
Thus $\overline V_{\omega}|_{P \ltimes V^\sharp} \simeq \Ind_{P \ltimes  ( V^+ \times (V^0)^\sharp)}^{P \ltimes V^\sharp} \left(\overline{V^0}_{\omega} \otimes ({\barFell}_{\chi^{V_1^+}} \ltimes 1)\right). $
\qed

Using this lemma, we can prove an intertwining result that is slightly stronger than the one needed to prove irreducibility in the complex setting. In order to formulate it, let $\KYu_{0+}=(G_1)_{x,\frac{r_1}{2}}(G_2)_{x,\frac{r_2}{2}} \hdots (G_n)_{x,\frac{r_n}{2}}(G_{n+1})_{x,0+}$.

\begin{Lemma}\label{Lemma-main}
	If $g \in G(F)$ such that $\Hom_{\KYu_{0+} \cap {^{g}\KYu}}(\kappa, {^g\repKYu}) \neq \{ 0 \}$, then $g \in \KYu$, where $^g\KYu$ denotes $g\KYu g^{-1}$ and $^g\repKYu(x)=\repKYu(g^{-1}xg)$.
\end{Lemma}
\Proof
The main part of the proof of Theorem 3.1 of \cite{Fi-Yu-works} consists of showing that in the setting with complex coefficients if $g \in G(F)$ such that 
\begin{equation*} \Hom_{\KYu \cap {^g\KYu}}\left( \repKYu|_{\KYu \cap {^g\KYu}}, {^g\repKYu|_{\KYu \cap {^g\KYu}}}\right) 
 \neq \{0 \} ,
\end{equation*}
then $g \in \KYu$. Using Lemma \ref{Lemma-Gerardin-mod-ell} instead of the equivalent result for complex coefficients, all the argument of the proof in the complex setting also work with $R$-coefficients to yield the same result in the mod-$\ell$ setting.

Now these arguments in turn can also be applied to show that if $g \in G(F)$ such that $\Hom_{\KYu_{0+} \cap {^{g}\KYu}}(\kappa, {^g\repKYu}) \neq \{ 0 \}$, then $g \in \KYu$. More precisely, suppose  
$$\Hom_{\KYu_{0+} \cap {^{g}\KYu}}(\kappa|_{\KYu_{0+} \cap {^{g}\KYu}}, {^g\repKYu}|_{\KYu_{0+} \cap {^{g}\KYu}}) \neq \{ 0 \}. $$
 The induction argument at the beginning of the proof of Theorem 3.1 of \cite{Fi-Yu-works} (which is essentially \cite[Corollary~4.5]{Yu}) only uses that  $\Hom_{\KYu_{0+} \cap {^{g}\KYu_{0+}}}(\repKYu, {^g\repKYu}) \neq \{ 0 \}$. Since $\repKYu|_{\KYu_{0+}}=(\kappa \otimes \rho)|_{\KYu_{0+}}$ is isomorphic to a direct sum of copies of $\kappa|_{\KYu_{0+}}$, these arguments show that $g \in \KYu G_{n+1}(F) \KYu $. 
 Since $\KYu$ normalizes $\KYu_{0+}$ and $\kappa$, we may therefore assume without loss of generality that $g \in G_{n+1}(F)$.

Our assumption that  $\Hom_{\KYu_{0+} \cap {^{g}\KYu}}(\kappa, {^g\repKYu}) \neq \{ 0 \}$ is equivalent to  $\Hom_{{\KYu} \cap ({^{g^{-1}}\KYu_{0+}})}({^{g^{-1}}\kappa}, {\repKYu}) \neq \{ 0 \}$, and $g \in G_{n+1}(F)$ implies $g^{-1} \in G_{n+1}(F)$. Now we can apply the arguments of the second half of the proof of Theorem 3.1 in \cite{Fi-Yu-works} to conclude that $g^{-1} \in (G_{n+1})_{[x]}$, hence $g \in (G_{n+1})_{[x]}$, as follows. The proof  uses a non-zero element $f \in \Hom_{{\KYu} \cap {^{g^{-1}}\KYu}}({^{g^{-1}}\repKYu}, {\repKYu})$ (where $g^{-1}$ is replaced by $g$ in the setting of \cite{Fi-Yu-works})
only to deduce that the image of $f$ is a non-trivial subspace on which ${^{g^{-1}}\KYu_{0+}}$ acts via a direct sum of copies of $^{g^{-1}}\kappa|_{{^{g^{-1}}\KYu_{0+}}}$. (More precisely, for the interested reader, we prove that the restriction of ${^{g^{-1}}\repKYu}$ (equivalently, the restriction of ${^{g^{-1}}\kappa}$) to a suitable subgroup $(H_{n+1})_{g^{-1}.x,0+}$ of ${^{g^{-1}}\KYu_{0+}}$ is given by $\phi=\prod_{i=1} \phi^i|_{(H_{n+1})_{g^{-1}.x,0+}}$ (times identity) and then use that therefore $(H_{n+1})_{g^{-1}.x,0+}$ acts via the character $\phi$ on the image of $f$.) 
 Hence the second half of the proof of Theorem 3.1 works equally well if we consider the image of a non-zero element in $\Hom_{{\KYu} \cap ({^{g^{-1}}\KYu_{0+}})}({^{g^{-1}}\kappa}, {\repKYu})$.

Thus we deduce that $\Hom_{\KYu_{0+} \cap {^{g}\KYu}}(\kappa, {^g\repKYu}) \neq \{ 0 \}$ implies $g \in \KYu$.
\qed

\textbf{Proof of Thereom \ref{Thm-mod-ell}.} \\
 We first note that $\repKYu$ is an irreducible representation of $\KYu$. The proof of this statement is exactly as in Lemma 3.3 of \cite{Fi-Yu-works}. By Lemma \ref{Lemma-main}, we have that
 if $g \in G(F)$ such that 
\begin{equation*}  \Hom_{\KYu \cap {^g\KYu}}\left( ^g\repKYu|_{\KYu \cap {^g\KYu}},  \repKYu|_{\KYu \cap {^g\KYu}}\right) \neq \{0 \} ,
\end{equation*}
then $g \in \KYu$. 
 Hence (by Mackey theory) we have $\End_G(\cind_{\KYu}^{G(F)}\repKYu)=R$. In order to show that $\cind_{\KYu}^{G(F)}\repKYu$  is irreducible, it suffices by Vigneras (\cite[Lemma~4.2]{Vigneras-irreducible}) to prove that for every smooth, irreducible $R$-representation $\pi'$ of $G(F)$ for which $\repKYu$ is contained in $\pi'|_{\KYu}$ the representation $\repKYu$ is also a quotient of $\pi'|_{\KYu}$. Let $\pi'$ be such a smooth, irreducible $R$-representation of $G(F)$ so that $\repKYu$ is contained in $\pi'|_{\KYu}$. 
  Since $\KYu_{0+}$ is a pro-$p$ group, we can write $\pi'|_{\KYu_{0+}} =\pi'_{\kappa} \oplus {\pi'}^\kappa$ where ${\pi'_\kappa}$ is $\kappa|_{\KYu_{0+}}$-isotypic, i.e. a direct sum of copies of the irreducible representation $\kappa|_{\KYu_{0+}}$, and none of the subquotients of ${\pi'}^\kappa$ contains $\kappa|_{\KYu_{0+}}$. As $\KYu_{0+}$ is a normal pro-$p$ subgroup of $\KYu$ and $\kappa$ is a representation of $\KYu$, the decomposition $\pi'|_{\KYu_{0+}}=\pi'_{\kappa} \oplus {\pi'}^\kappa$ is preserved by $\KYu$, i.e. we have a decomposition  $\pi'|_{\KYu} = \pi'_{\kappa} \oplus {\pi'}^\kappa$ of $\KYu$-representations. Similarly, we have a decomposition 
\begin{equation} \label{eqn-decomposition} \bigoplus_{ g  \in \, \KYu\backslash G(F)/\KYu} \ind_{\KYu \cap {^g\KYu}}^{\KYu}({^g\repKYu}|_{\KYu \cap {^g\KYu}}) \simeq \pi|_{\KYu} = \pi_{\kappa} \oplus \pi^\kappa .
\end{equation}
By Mackey theory, 
$$  \left(\ind_{\KYu \cap {^g\KYu}}^{\KYu}{^g\repKYu}\right)|_{\KYu_{0+}}= \bigoplus_{ k \in \,  \KYu_{0+}\backslash \KYu/(\KYu \cap {^g\KYu})} \ind^{\KYu_{0+}}_{\KYu_{0+} \cap {^k(\KYu \cap {^g\KYu})}} {^{kg}\repKYu} = \bigoplus_{ k \in \,  \KYu_{0+}\backslash \KYu/(\KYu \cap {^g\KYu})} \ind^{\KYu_{0+}}_{\KYu_{0+} \cap {^{kg}\KYu}} {^{kg}\repKYu} .$$
By Frobenius reciprocity, $\ind^{\KYu_{0+}}_{\KYu_{0+} \cap {^{kg}\KYu}} {^{kg}\repKYu}$ contains $\kappa|_{\KYu_{0+}}$ if and only if $\Hom_{\KYu_{0+} \cap {^{kg}\KYu}}(\kappa,{^{kg}\repKYu})\neq \{0\}$.
By  Lemma \ref{Lemma-main} this can only happen if $kg \in \KYu$, hence $g \in \KYu$. Thus Equation \eqref{eqn-decomposition} yields $\pi_\kappa \simeq \repKYu$, because $\repKYu = \rho \otimes \kappa$. Since $\repKYu$ is contained in $\pi'|_{\KYu}$ and $\pi'$ is irreducible, we obtain by Frobenius reciprocity a surjection from $\pi=\cind_{\KYu}^{G(F)}\repKYu$ onto $\pi'$ that maps $\pi_\kappa \simeq \repKYu$ surjectively onto $\pi'_{\kappa}$. Recall that we assumed that $\pi'|_{\KYu}$ contains the irreducible representation $\repKYu$. Hence $\pi'_{\kappa} \simeq \repKYu$, and therefore $\repKYu$ is a quotient of $\pi'|_{\KYu}$.
Thus $\cind_{\KYu}^{G(F)}\repKYu$ is irreducible. As the matrix coefficients of $\cind_{\KYu}^{G(F)}\repKYu$ are compactly supported mod center, we also obtain from \cite[II.2.7]{Vigneras-book} that $\cind_{\KYu}^{G(F)}\repKYu$ is cuspidal. \qed

\section{Exhaustion of cuspidal representations} \label{Section-exhaustion}

\begin{Thm} \label{Thm-exhaustion-mod-ell}
	Assume that $p$ does not divide the order of the (absolute) Weyl group of $G$. Then every smooth, irreducible, cuspidal $R$-representation of $G(F)$ arises from the construction in Section \ref{Section-construction-mod-ell}.
\end{Thm}

\begin{Rem}
	This theorem was proven for representations with complex coefficients in \cite{Fi-exhaustion}. In \cite{Fi-exhaustion} we proved the more general result that every irreducible, smooth complex representation contains an $\fs$-type (for some inertial equivalence class $\fs$). Recall that Bernstein (\cite{Bernstein}) introduced a decomposition of the category of smooth complex representations of $G(F)$ into blocks indexed by the cuspidal support, and these blocks can also be characterized by $\fs$-types (see \cite{BK-types}). However, for $R$-representations there are situations where the cuspidal support does not produce a decomposition of the category of smooth $R$-representation into blocks. We can nevertheless use the approach from \cite{Fi-exhaustion} applied to cuspidal $R$-representations to prove an exhaustion theorem for irreducible cuspidal $R$-representations by adding a few additional arguments.
\end{Rem}

Before we proof the theorem, let us introduce some notation. Let $G=G_1 \supseteq G_2 \supsetneq G_3 \supsetneq \hdots \supsetneq G_{n+1}$ be a sequence of twisted Levi subgroups of $G$ that split over a tamely ramified extension of $F$ and $x$ a point in $\sB(G_{n+1},F) \subset \sB(G,F)$. Then Moy and Prasad (\cite[6.3 and 6.4]{MP2}) attach to $x$ and $G_{n+1}$ a Levi subgroup $M_{n+1}$ of $G_{n+1}$ such that $x \in \sB(M_{n+1},F) \subset \sB(G_{n+1},F)$ and $(M_{n+1})_{x,0}$ is a maximal parahoric subgroup of $M_{n+1}(F)$ with $(M_{n+1})_{x,0}/(M_{n+1})_{x,0+} \simeq (G_{n+1})_{x,0}/(G_{n+1})_{x,0+}$. Following \cite[2.4]{Kim-Yu} we denote by $Z_s(M_{n+1})$ the maximal split torus in the center of $M_{n+1}$ and by $M_i$ the centralizer of $Z_s(M_{n+1})$ in $G_i$ for $1 \leq i \leq n$.

\textbf{Proof of Theorem \ref{Thm-exhaustion-mod-ell}.}\\
Let $(\pi', V')$ be a smooth, irreducible, cuspidal $R$-representation of $G(F)$. We like to apply a version of \cite[Theorem~7.12]{Fi-exhaustion} in our mod-$\ell$ setting. In order to do so, recall that in Section 7 of \cite{Fi-exhaustion}, we use Pontryagin duality to show that a smooth, complex unitary character of a compact open subgroup $A_{0+}$ of a locally compact abelian group $A$ extends to a unitary complex character of the whole group $A$ (in the notation of \cite{Fi-exhaustion} we use it for the pair $A=(G_{j+1}/H_{j+1})(F)$ and $A_{0+}=(G_{j+1}/H_{j+1})(F)_{0+}$ and for the pair  $A=(G_{j}/H_{j})(F)$ and $A_{0+}=(G_{j}/H_{j})(F)_{r_j+}$). In our setting, the smooth character of $A_{0+}$ takes values in the roots of unity of $\barFell^* \subset R^*$. Since the group of roots of unity in $\barFell^*$ is divisible, our character of $A_{0+}$ extends to an $\barFell$-character of $A$.

Now we can use all the arguments of \cite[Theorem~7.12]{Fi-exhaustion} 
 that also apply in the mod-$\ell$ setting to obtain from  $(\pi', V')$ a tuple $((G_i)_{1 \leq i \leq n+1}, x, (r_i)_{1 \leq i \leq n}, \rho^M, (\phi_i)_{1 \leq i \leq n})$ where
\begin{enumerate}[label=(\alph*),ref=\alph*]
	\item $G=G_1 \supseteq G_2 \supsetneq G_3 \supsetneq \hdots \supsetneq G_{n+1}$ are twisted Levi subgroups of $G$ that split over a tamely ramified extension of $F$
	\item  $x \in \sB(G_{n+1},F)\subset \sB(G,F)$ 
	\item $r_1 > r_2 > \hdots > r_n >0$ are real numbers
	\item[\refstepcounter{enumi}(\alph{enumi}')] $\rho^M$ is an irreducible $R$-representation of $(M_{n+1})_{x}$ that is trivial on $(M_{n+1})_{x,0+}$
	\item $\phi_i$, for $1 \leq i \leq n$, is a character of $G_{i+1}(F)$ of depth $r_i$ that is trivial on $(G_{i+1})_{x,r_i+}$ 
\end{enumerate}
satisfying the following conditions 
\begin{enumerate}[label=(\roman*),ref=\roman*]
	\refstepcounter{enumi}
	\item[\refstepcounter{enumi}(\roman{enumi}')] the point $x \in \sB(M_{n+1},F) \subset \sB(G_{n+1},F)\subset \sB(G,F)$ satisfies  $$ \sum_{j=1}^{n} \left(\dim((G_i)_{x,\frac{r_i}{2}}/(G_i)_{x,\frac{r_i}{2}+}) - \dim((M_i)_{x,\frac{r_i}{2}}/(M_i)_{x,\frac{r_i}{2}+}) \right)=0  $$ 
	\item[\refstepcounter{enumi}(\roman{enumi}')] $\rho^M|_{(M_{n+1})_{x,0}}$ is a cuspidal representation of $(M_{n+1})_{x,0}/(M_{n+1})_{x,0+} \simeq (G_{n+1})_{x,0}/(G_{n+1})_{x,0+}$
	\item $\phi_i$ is $G_i$-generic of depth $r_i$ relative to $x$ for all $1 \leq i \leq n$ with $G_i \neq G_{i+1}$
\end{enumerate}
with the following property:

Following Kim and Yu (\cite[7.1 and 7.3]{Kim-Yu}) we define the group $K_{G_{n+1}}$ to be the group generated by $(M_{n+1})_x$ and $(G_{n+1})_{x,0}$ and let $K$ be the compact open subgroup 
$$K =(G_1)_{x,\frac{r_1}{2}}(G_2)_{x,\frac{r_2}{2}} \hdots (G_n)_{x,\frac{r_n}{2}}K_{G_{n+1}} 
\subset G(F).$$
Then the construction in Section \ref{Section-construction-mod-ell} applied to $K$ instead of $\KYu$ yields an irreducible representation $\repKYu_K  =\rho^M_K \otimes \kappa_K$ of $K$, and we have $\Hom_K(\repKYu_K,\pi'|_K) \neq \{0\}$. (In the case of complex coefficients the pair $(K, \repKYu_K)$ would be a type by \cite[7.5~Theorem and 7.3.~Remark]{Kim-Yu}.) 

By construction, the restriction of $\kappa_K$ to the subgroup 
$$K_+= (G_1)_{x,\frac{r_1}{2}+}(G_2)_{x,\frac{r_2}{2}+} \hdots (G_n)_{x,\frac{r_n}{2}+}(G_{n+1})_{x,0+}$$
is given by a character $\hat \phi$ times the identity.
 In the proof of \cite[6.3.~Theorem]{Kim-Yu} Kim and Yu show that the Jacquet functor $r_{{M_1}, G}: V' \ra V'_{M_1}$ induces an injection on ${V'}^{(K_+, \hat \phi)}$, where ${V'}^{(K_+, \hat \phi)}$ denotes the subspace of $V'$ on which $K_+$ acts via the character $\hat\phi$.
 While Kim and Yu work with complex coefficients, their proof also works in our setting. More precisely, first note that while our point $x$ is a point of $\sB(G,F)$ that lies on the image of (any) inclusion of $\sB(M_{n+1},F)$ into $\sB(G,F)$, Kim and Yu start with a point $x$ in $\sB(M_{n+1},F)$ and include a diagram of embeddings of $\sB(M_i, F)$ and $\sB(G_i, F)$ into $\sB(G, F)$ ($1 \leq i \leq n+1$) as part of their input datum. Kim and Yu's requirement that their diagram of embeddings of buildings is, using their notation, ``$\vec s$-generic relative to $x$'' corresponds to our Condition (ii').
  Now the proof of Kim and Yu follows the strategy of the proof by Moy and Prasad (\cite{MP2}) in the depth-zero case, i.e. the special case that $n=0$, $K_+=G_{x,0+}$ and $\hat \phi=1$, and also relies on Moy and Prasad's result (part of \cite[Proposition~6.7]{MP2}) as induction hypothesis. However, the proof of Moy and Prasad works mod $\ell$. It in turn relies on a result of Howlett and Lehrer (\cite{Howlett-Lehrer}) who cover the mod-$\ell$ case.

Since $\repKYu_K|_{K_+} =(\rho^M_K)|_{K_+} \otimes (\kappa_K)|_{K_+} =\Id \otimes (\hat \phi \cdot \Id)$ and $\Hom_K(\repKYu_K,\pi'|_K) \neq \{0\}$, the subspace ${V'}^{(K_+, \hat \phi)}$ is nonzero. Hence the image of the Jacquet functor $r_{{M_1}, G}V'$ is nonzero, and therefore, since $(\pi', V')$ is cuspidal, we obtain that ${M_1}=G$. Moreover, this implies that $Z_s(M_{n+1}) \subset Z(G)$, and hence $Z(G_{n+1})/Z(G)$ is anisotropic and $M_{n+1}=G_{n+1}$, which implies that the image of $x$ in $\sB(G_{n+1}^{\mathrm{der}},F)$ is a vertex. 
Therefore, by working with $(G_{n+1})_{[x]}$ instead of $K_{G_{n+1}}=(G_{n+1})_{x}$ in the proof of \cite[Theorem~7.12]{Fi-exhaustion}, we obtain a tuple $((G_i)_{1 \leq i \leq n+1}, x, (r_i)_{1 \leq i \leq n}, \rho, (\phi_i)_{1 \leq i \leq n})$ as in Section \ref{Section-input-mod-ell} that allows us to construct a smooth, irreducible, cuspidal $R$-representation $\pi=\cind_{\KYu}^{G(F)}\repKYu$ of $G(F)$ such that $\repKYu$ is contained in $\pi'|_{\KYu}$. By Frobenius reciprocity we obtain a non-trivial morphism between the two irreducible representations $\pi$ and $\pi'$, and hence $\pi'\simeq \cind_{\KYu}^{G(F)}\repKYu$. \qed

\bibliography{Fintzenbib}  

\end{document}